\RequirePackage{fix-cm}
\documentclass[11pt]{amsart}
\usepackage[margin=4cm]{geometry}
\usepackage{amsmath,amsfonts,amsthm,amssymb,amscd,verbatim}
\usepackage{graphicx,color,booktabs,bm,microtype}
\usepackage{colortbl,makecell}
\usepackage[pagebackref]{hyperref}
\numberwithin{equation}{section}
\newtheorem{theorem}{Theorem}[section]
\newtheorem{lemma}[theorem]{Lemma}

\newtheorem{proposition}[theorem]{Proposition}
\newtheorem{corollary}[theorem]{Corollary}

\newtheorem{conjecture}{Conjecture}

\theoremstyle{definition}

\theoremstyle{remark}
\newtheorem{remark}[theorem]{Remark}

\newcommand{\Alt}{\mathop{\mathrm{Alt}}}

\newcommand{\Aut}{\mathop{\mathrm{Aut}}}

\newcommand{\GL}{\mathop{\mathrm{GL}}}

\newcommand{\F}{\mathbb F}
\newcommand{\Z}{\mathbb Z}

\newcommand{\E}{\mathbb E}
\newcommand{\one}{\mathbf 1}
\newcommand{\abs}[1]{\lvert#1\rvert}

\DeclareMathOperator{\Soc}{Soc}
\newcommand{\arxiv}[1]{\href{https://arxiv.org/abs/#1}{arXiv:#1}}

\definecolor{mycolor}{rgb}{0.55,0.0,0.16}
\definecolor{myred}{rgb}{0.75,0.0,0.16}
\definecolor{mygreen}{rgb}{0.0,0.4,0.16}
\definecolor{myviolet}{rgb}{1,0,1}
\definecolor{mypink}{rgb}{0.67,0,0.47}

\hypersetup{colorlinks=true,linktocpage=true,pageanchor=true,
  linkcolor=mycolor,urlcolor=mypink,citecolor=mygreen,
  pdftitle={Progression-free sets in class 2 and arbitrary finite groups},
  pdfauthor={Luca Sabatini}}

\author[Luca Sabatini]{Luca Sabatini}
\address{Luca Sabatini, Mathematics Institute, University of Warwick}
\email{luca.sabatini@warwick.ac.uk, sabatini.math@gmail.com}
\title[Progression-free sets in arbitrary finite groups]{Progression-free sets in
class $2$\\and
arbitrary finite groups}

\begin{document}

\begin{abstract}
We study progression-free sets in arbitrary finite groups.
We transfer sharp upper bounds from abelian groups to nilpotent groups of class $2$,
and from nilpotent groups to arbitrary groups.
Combining a recent result of Raghavan with the author's theorem on large subgroups of class $2$ gives
\[
 r_3(G) \, \leq \, C \, |G|\exp\left(
 -c\frac{(\log|G|)^{1/6}}{(\log\log|G|)^{1/3}}\right)
\]
for every finite group as $|G|\to\infty$.
This improves the estimates obtained by using the abelian subgroup supplied by Pyber's theorem.
\end{abstract}

\maketitle


\section{Introduction}\label{sec1}

One of the earliest works in additive combinatorics is a paper of Erd\H{o}s and Tur\'an \cite{ET36} from 1936. They defined $r_k(N)$ as the cardinality of the largest subset of
$\{1,\ldots,N\}$ not containing $k$ distinct elements in arithmetic progression, and conjectured that $r_k(N)=o(N)$ for any fixed $k \geq 3$.
For $k=3$, this was proved by Roth \cite{Rot53},
and Sanders \cite{San11} obtained
$r_3(N) \ll N/(\log N)^{1-o(1)}$ more than half a century later.
The breakthrough of Kelley and Meka \cite{KM23}, and the subsequent work of Bloom
and Sisask \cite{BS23}, gave bounds with decay
$\exp(-c(\log N)^\alpha)$ with $\alpha=1/9$.
The current best estimate is due to Raghavan \cite{Rag26}, who proved
\[ r_3(N) \> \ll \> N \exp\left(
 -c\left(\frac{\log N}{\log\log N}\right)^{1/6}\right) . \]
In the opposite direction, Behrend's celebrated construction \cite{Beh46} (see also \cite{EHPS24}) gives
\[ r_3(N) \> \ge \> N \exp(-c\sqrt{\log N}) , \]
and it seems a common belief that this is the right order of magnitude for $r_3(N)$.

For fixed $k\geq4$, Szemer\'edi \cite{Sze75} proved $r_k(N)=o(N)$ in 1975.
Gowers \cite{Gow01} introduced higher-order Fourier analysis and obtained
$r_k(N)\ll_k N/(\log\log N)^{c_k}$ for some $c_k>0$.
In a preprint dated 23 September 2026, OpenAI \cite{OAI26} claimed a proof that there are constants $c_k,\alpha_k>0$ such that
\[
 r_k(N) \> \ll_k \> N\exp\bigl(-c_k(\log N)^{\alpha_k}\bigr).
\]
This matches Rankin's extension of Behrend's construction \cite{Ran61}.
For a more detailed historical survey of arithmetic progressions, see \cite{Gre26}.

Bounds over the integers are closely related to those for finite cyclic groups.
Taking representatives in $\{1,\ldots,N\}$ of a subset of $\Z/N\Z$ preserves the property of avoiding progressions,
and a $k$-progression-free subset of $\{1,\ldots,N\}$ remains $k$-progression-free when embedded in $\Z/kN\Z$.
On the other hand, torsion can play an important role for noncyclic abelian groups.
For example, Ellenberg and Gijswijt \cite{EG17} proved that, for every
fixed odd prime $p$, every $3$-progression-free subset of $\F_p^d$ has size
$O_p(c_p^d)$ for some $c_p<p$.

In this paper we consider the non-commutative setting.
A {\em $k$-progression} in an arbitrary group $G$ is
\[
 (x,xy,\ldots,xy^{k-1}),\qquad x,y\in G ,
\]
and it is \emph{nontrivial} if $y\ne1$. A subset of $G$ is
\emph{$k$-progression-free} if it contains no nontrivial $k$-progression, and we set
\[
 r_k(G) \> := \> \max\{|X|:X\subseteq G\text{ is $k$-progression-free}\} .
\]

Clearly $r_k(G) \geq r_3(G)$ for any $k \geq 3$.
Non-commutative Roth-type problems were first investigated by Bergelson, McCutcheon and Zhang \cite{BMZ97}.
Solymosi proved that subsets of fixed positive density contain a progression with three distinct terms whenever the odd part of $|G|$ is sufficiently large \cite[Theorem~2.5]{Sol13}.
In full generality, Sanders \cite[Theorem~1.4]{San17} applied a finite-abelian bound in a coset of the abelian subgroup supplied by Pyber \cite{Pyb97}.
This strategy comes at the cost of halving the exponent \(\alpha\) in the density term \(\exp(-c(\log |G|)^\alpha)\).

The principal novelty in our approach is to use substantially more group theory.
With the conviction that the intervals should represent the most difficult case,
we propose the following conjecture for arbitrary finite groups.

\begin{conjecture} \label{conj:A}
For every $k \geq 3$ there is $\theta_k >0$ such that every finite group $G$ satisfies
\[
 \frac{r_k(G)}{|G|} \> \ll_k \> \left(\frac{r_k(|G|)}{|G|}\right)^{\theta_k} .
\]
\end{conjecture}

We will come close to prove Conjecture~\ref{conj:A}, especially for $k=3$, which we consider our main contribution. However, for technical reasons, we prefer to state most of our results in a different form.
Given $0<\alpha\leq1/2$, we say that a class $\mathcal{C}$ of finite
groups satisfies {\em property~$(\mathrm{AP}_{k,\alpha})$} if there are constants $c,C>0$,
such that for every $G\in\mathcal{C}$ we have
\[
 r_k(G) \> \leq \>  C|G|\exp\bigl(-c(\log|G|)^\alpha\bigr).
\]

The starting point of our investigation is to replace the abelian subgroup of Pyber with the larger nilpotent subgroup of class $2$ supplied by \cite{Sab22}.
Along these lines, our first main result transfers bounds from abelian groups to groups of class $2$.

\begin{theorem}[From abelian to class $2$] \label{thm:abeTOclass2}
Fix $k \geq 3$.
If the class of finite abelian groups satisfies $(\mathrm{AP}_{k,\alpha})$ for some
$0<\alpha\leq1/2$, then so does the class of finite nilpotent groups
of class at most $2$.
\end{theorem}

Given a group $H$ of class $2$, we construct a commutative group law on the underlying set of a subgroup $H_0$ of $H$ that preserves progressions.
For groups of odd order we have $H=H_0$ and in fact this is the classical trick of Baer \cite{Bae38}.

With Theorem~\ref{thm:abeTOclass2} in hand, \cite{Sab22} allows us to pass from abelian to arbitrary finite groups with only minor quantitative loss.
If we instead assume $(\mathrm{AP}_{k,\alpha})$ for nilpotent groups of arbitrary class, we can further remove this loss and obtain the following smooth statement.

\begin{theorem}[From nilpotent to arbitrary groups] \label{thm:nilpTOarbitrary}
Fix $k \geq 3$.
If the class of finite nilpotent groups satisfies $(\mathrm{AP}_{k,\alpha})$ for some
$0<\alpha\leq1/2$, then so does the class of all finite groups.
In particular, Conjecture~\ref{conj:A} holds if and only if it holds
for finite nilpotent groups.
\end{theorem}

In the proof of Theorem~\ref{thm:nilpTOarbitrary},
we have to deal with alternating composition factors in an efficient way.
We use a result of Gonz\'alez-S\'anchez and Klopsch \cite{GSK11} to create a large $3$-subgroup from a large elementary abelian section,
and then apply the cap-set bound \cite{EG17} to this subgroup.

So far, all our results are conditional on bounds in the abelian case.
For $k=3$, we leverage \cite{EG17} to further reduce the problem to an estimate over the intervals (Proposition~\ref{prop:intervTOabe}).
This allows us to apply the theorem of Raghavan \cite{Rag26}
and to obtain the following close bound for all finite groups.

\begin{theorem} \label{thm:unconditional}
There is an absolute constant $c>0$ such that every finite group $G$ satisfies
\[
 r_3(G) \> \ll \> |G|\exp\left(
 -c\frac{(\log|G|)^{1/6}}{(\log\log|G|)^{1/3}}\right).
\]
In particular, the class of all finite groups satisfies $(\mathrm{AP}_{3,\alpha})$
for every $0<\alpha<1/6$.
\end{theorem}

\begin{table}[htbp]
\centering
\small
\renewcommand{\arraystretch}{1.8}
\setlength{\tabcolsep}{5pt}
\setlength{\arrayrulewidth}{0.8pt}
\begin{tabular}{|>{\columncolor[gray]{0.9}}c|c|c|}
\hline
\rowcolor[gray]{0.75}
\cellcolor[gray]{0.5} & $k=3$ & $k\geq4$ \\
\hline
Intervals $\to$ Abelian groups
& Proposition~\ref{prop:intervTOabe}
& Open \\
\hline
Abelian groups $\to$ Class $2$
& \multicolumn{2}{c|}{Theorem~\ref{thm:abeTOclass2}} \\
\hline
Class $2$ $\to$ Nilpotent groups
& \multicolumn{2}{c|}{%
  \makecell[c]{minor loss;\\[3pt]
               removing the loss is open}} \\
\hline
Nilpotent groups $\to$ Arbitrary groups
& \multicolumn{2}{c|}{Theorem~\ref{thm:nilpTOarbitrary}} \\
\hline
\end{tabular}
\vspace{0.2cm}
\caption{Partial results towards Conjecture~\ref{conj:A}.}
\label{tab:transfers}
\end{table}

The paper is organized as follows.
In Section \ref{sec:preliminaries} we give the necessary background, and in Section~\ref{sec:abeTOclass2} we prove Theorem~\ref{thm:abeTOclass2}. In Section~\ref{sec:arbitrary} we move to arbitrary finite groups and prove Theorem~\ref{thm:nilpTOarbitrary}.
In Section~\ref{sec:unconditional}, we restrict to $k=3$ and prove~Theorem~\ref{thm:unconditional}.

\vspace{0.2cm} \noindent
{\bfseries Notation.} Unsubscripted logarithms are natural. 
Given two positive functions $f$ and $g$, by $f \ll g$ and $f= O(g)$ we
mean the same thing, namely that there is a constant $C$ such that $f \leq Cg$.
Implicit constants are absolute unless a dependence is specified.

\vspace{0.2cm}\noindent
{\bfseries Acknowledgments.}
Part of this work was undertaken during a visit to the mathematics department of the University of Cambridge.
The author would like to thank P\'eter Varj\'u for his hospitality, advice, and constant availability to discuss mathematics.

\vspace{0.2cm}\noindent
{\bfseries AI statement.}
The author's first encounter with Lemma~\ref{lem:coset} was while reading Sanders' paper \cite{San17}, and his familiarity with large-subgroup results such as Theorem~\ref{thm:large-subgroups} motivated the present work.
The idea of using the Baer trick on a subgroup was discovered during an initial conversation with an LLM (actually the original argument was quite different and
artificial;
in particular, Baer was not mentioned at all).
None of the other mathematical ideas in this paper is generated by an LLM.

\vspace{0.1cm}
\section{Preliminaries} \label{sec:preliminaries}

\subsection{Progression-free sets and subgroups}

The projection of a progression-free set in a quotient need not be progression-free.
On the other hand, progressions are preserved by left translation.
This allows passage to subgroup cosets, as exploited in the following easy but important observation (see \cite[Proof of Theorem~1.4]{San17}).
For convenience,
let us write
\[ \rho_k(G) \> := \> \frac{r_k(G)}{|G|} \]
for the maximum density of a $k$-progression-free set.

\begin{lemma}[Subgroup density inequality] \label{lem:coset}
If $H \leqslant G$, then $\rho_k(G) \leq \rho_k(H)$.
\end{lemma}

\begin{proof}
Let $X \subseteq G$ be $k$-progression-free.
The average of $|X\cap gH|$ over the distinct left cosets of $H$ in $G$ is
$|X|/|G:H|$. In particular, there exists $g \in G$ such that $|X \cap gH| \geq |X|/|G:H|$. The claim follows because $g^{-1}X \cap H$ is $k$-progression-free in $H$.
\end{proof}

A substantial portion of our work in this paper consists in applying Lemma~\ref{lem:coset} to convenient subgroups.
The most extreme example is the following.

\begin{corollary} \label{cor:involutions}
If $E\leqslant G$ is an elementary
abelian $2$-subgroup, then
$\rho_k(G) \leq |E|^{-1}$.
\end{corollary}

\begin{proof}
Note that $r_k(E)=1$, because any two distinct elements $a,b \in E$ form an arbitrarily long progression $(a,b,a,b,\ldots)$.
\end{proof}

 We now record a variant of the estimate in \cite{EG17},
 in which the density bound is uniform in the prime.

\begin{theorem}[Ellenberg--Gijswijt] \label{thm:EG}
There is an absolute constant $c>0$ such that
for every prime $p$ and every integer $d\geq1$ we have
\[ \rho_3(\F_p^d) \> \ll \> e^{-c \, d}. \]
\end{theorem}

\begin{proof}
For $p=2$, the claim is trivial because $r_3(\F_2^d)=1$.

For $p$ odd, apply \cite[Theorem~4]{EG17} with coefficients $(1,-2,1)$, so that the equation is $x+z=2y$.
We obtain
\[ r_3(\F_p^d) \leq 
3 \> \left| \left\{ (a_1,\ldots,a_d) : 
a_i \in \{0,\ldots,p-1\} , \sum_{i=1}^d a_i \leq \frac{d(p-1)}{3} \right\} \right| . \]
This can be rewritten as
\[
 \rho_3(\F_p^d)
 \leq 3 \,\mathbb P\left(
       U_1+\cdots+U_d\leq\frac{d(p-1)}3\right) ,
\]
where the $U_i$'s are independent random variables uniform on $\{0,\ldots,p-1\}$.
Hoeffding's inequality \cite{Hoe63} gives
\[
 \mathbb P\left(U_1+\cdots+U_d\leq\frac{d(p-1)}3\right)
 \leq \exp\left(-\frac{2(d(p-1)/6)^2}{d(p-1)^2}\right)
 \leq e^{-c \, d}
\]
for an absolute constant $c>0$.
\end{proof}

We are not aware of similarly strong bounds when $k \geq 4$.

\subsection{Large subgroups}

A finite group $H$ is nilpotent if it is the direct product of finite $p$-groups.
It has nilpotency class at most $2$ if every commutator is central, i.e. if $[H,H] \subseteq Z(H)$. Every abelian group satisfies this property trivially.

The following bounds are the best possible and are valid in every finite group.

\begin{theorem}[Pyber \cite{Pyb97}; Sabatini \cite{Sab22}]
\label{thm:large-subgroups}
There is an absolute constant $c_0>0$ such that every finite group
$G$ contains an abelian subgroup $A$ and a nilpotent subgroup $H$ of
class at most $2$ satisfying
\[ |A| \, \ge \, \exp\bigl(c_0\sqrt{\log|G|}\bigr),\qquad
 |H| \, \ge \, \exp\left(c_0\frac{\log|G|}{\log\log|G|}\right) . \]
\end{theorem}

In fact the subgroup $H$ provided by \cite{Sab22} is \emph{ab-maximal}; such groups are studied further in \cite{ES25}.
For other recent applications of Pyber's subgroup to additive combinatorics, see \cite{Sch25}.

\begin{remark} \label{rem:Olsh}
A surprising non-explicit construction of Ol'shanskii \cite{Ols78} gives an arbitrarily large finite group $G$, which is also nilpotent of class $2$, without abelian subgroups larger than $\exp\left(\sqrt{8 \log|G|}\right)$.
\end{remark}

\subsection{Baer trick} 

Let us recall briefly the classical construction of Baer \cite[Appendix~B]{Bae38},
which is also explained in \cite[Lemma~4.37]{Isa08}.

Let $H$ be a nilpotent group of class $2$, and let $a,b \in H$.
We use the convention $[a,b]=a^{-1}b^{-1}ab$, so that
 \[  ab \> = \> ba [a,b] . \]
The error term $[a,b]$ is central and can be moved around.
If there is an element $\beta \in H$ such that $\beta^2 = [a,b]$, then it is possible to define
\[ a \oplus b \> = \> ab \beta^{-1} , \]
so that $a \oplus b = ab \beta^{-1} = ba \beta = b \oplus a$.
The existence of such an element $\beta$ is ensured if $H$ has odd order,
and the main challenge in the next section is to restrict the operation $\oplus$ to a meaningful subgroup $H_0$ of $H$ where commutators are squares. If $H$ has odd order, we will in fact have $H=H_0$.

\vspace{0.1cm}
\section{A Baer trick in even order}
\label{sec:abeTOclass2}

In this section we prove Theorem~\ref{thm:abeTOclass2}.
The key step is the following.

\begin{proposition} \label{prop:abelianCreation}
Fix $k\geq3$. Let $H$ be a finite nilpotent group of class $2$ and let $s= -\log_2(\rho_k(H))$.
Then there exists a finite abelian group $A$ such that
\[
  |A| \, \ge \, |H| \,4^{-s^2} ,
 \qquad
 \rho_k(H) \, \leq \, \rho_k(A) .
\]
\end{proposition}

Let $H^n$ be the subgroup generated by the $n$th powers.
The underlying set of the group $A$ in Proposition~\ref{prop:abelianCreation} is the subgroup $H^2Z(H)$ of $H$, and the density bound is just Lemma~\ref{lem:coset} in disguise.
As expected, the group law on $A$ is a generalization of the Baer trick and differs from multiplication in $H$
(also, note that such a large abelian subgroup would generally be impossible by Remark~\ref{rem:Olsh}).

\begin{lemma} \label{lem:beta}
Let $H$ be nilpotent of class $2$, $Z=Z(H)$, and let $H_0=H^2Z$.
Let $a,b \in H_0$, so that $a=u^2z$ and $b=v^2t$ for some $u,v\in H$ and $z,t\in Z$.
Define
\begin{equation} \label{eq:beta}
 \beta(a,b) \> := \> [u,v]^2 .
\end{equation}
Then $\beta \colon H_0 \times H_0 \to [H,H]$ is well-defined and
bimultiplicative. Moreover,
\[
 \beta(b,a)=\beta(a,b)^{-1},\qquad
 \beta(a,b)^2=[a,b].
\]
\end{lemma}

\begin{proof}
For well-definedness, suppose $u'$ is another representative with
$(u')^2Z=u^2Z$, say $u'=uw$. Then $w^2\in Z$, and the class-$2$
commutator identities give
\[
 [uw,v]^2=[u,v]^2[w,v]^2
          =[u,v]^2[w^2,v]=[u,v]^2.
\]
Changing $v$ is handled in the same way, and central changes in $z,t$
have no effect.

If $a_i=u_i^2z_i$ for $i=1,2$, then $a_1a_2$ has square root
$u_1u_2$ modulo $Z$. Hence
\[
 \beta(a_1a_2,b)=[u_1u_2,v]^2
   =[u_1,v]^2[u_2,v]^2
   =\beta(a_1,b)\beta(a_2,b).
\]
The other variable is analogous. When the first argument is in $Z$,
its square root $u$ modulo $Z$ satisfies $u^2\in Z$, so
$[u,v]^2=[u^2,v]=1$; the second argument is similar.
Choosing the same square root for both arguments proves
$\beta(a,a)=1$, and skew-symmetry follows. Finally,
\[
 [a,b]=[u^2,v^2]=[u,v]^4=\beta(a,b)^2
\]
as desired.
\end{proof}

\begin{lemma}[Generalized Baer trick] \label{lem:law}
Let $H$ be nilpotent of class $2$, $Z=Z(H)$, and let $H_0=H^2Z$.
Let $\beta \colon H_0 \times H_0 \to [H,H]$ be as in \eqref{eq:beta}, and define
\[
 a\oplus b \> = \> ab\,\beta(a,b)^{-1}\qquad(a,b\in H_0).
\]
Then $(H_0,\oplus)$ is an abelian group with identity $1$, inverse
$\ominus a=a^{-1}$, and integer multiples equal to powers in $H_0$.
Moreover, the $k$-progressions in $H_0$ and $(H_0,\oplus)$
coincide for every $k\geq3$. In particular, if $H_0$ is finite, then
\[ r_k(H_0) \> = \> r_k(H_0,\oplus) \]
for every $k \geq 3$.
\end{lemma}

\begin{proof}
Every nontrivial element of $H/H^2$ is an involution, so $H/H^2$ is abelian and $[H,H] \subseteq H^2$. Closure follows.
By Lemma~\ref{lem:beta} we have
\[
 ab\,\beta(a,b)^{-1}
  =ba[a,b]\beta(a,b)^{-1}
  =ba\,\beta(a,b)
  =ba\,\beta(b,a)^{-1},
\]
so $\oplus$ is commutative.
Since $\beta$ is bimultiplicative and trivial on central arguments,
both $(a\oplus b)\oplus c$ and $a\oplus(b\oplus c)$ are equal to
\[
 abc\,\beta(a,b)^{-1}\beta(a,c)^{-1}\beta(b,c)^{-1}.
\]
The assertions about the identity, inverses, and integer multiples can be checked directly.
In general, an ordered tuple $(a_0,\ldots,a_{k-1})$ is a $k$-progression
if and only if the successive ratios $a_i^{-1}a_{i+1}$ are equal, i.e.
\[
 a_{i+2}=a_{i+1}a_i^{-1}a_{i+1}
 \qquad(0\leq i\leq k-3).
\]

In our group $A$ we have
\begin{align*}
 (b\oplus a^{-1})\oplus b
 &=ba^{-1}b\,\beta(b,a^{-1})^{-1} \,\beta(ba^{-1},b)^{-1}\\
 &=ba^{-1}b\,\beta(a^{-1},b) \,\beta(a^{-1},b)^{-1}\\
 &=ba^{-1}b.
\end{align*}
Thus, $c=ba^{-1}b$ in $H_0$ is equivalent to $c= b \oplus a^{-1} \oplus b$ in $A$, as desired.
\end{proof}

\begin{remark}
The same construction works by choosing $H_0 = H^2$, but we prefer to use $H^2Z(H)$ because it is more natural and much bigger in general.
\end{remark}

We proceed bounding the index of $H^2$ in $H$ when $\rho_k(H)$ is close to $1$.
Let $d(G)$ be the minimum number of generators.

\begin{lemma} \label{lem:index}
Let $H$ be a finite group and let $s= -\log_2(\rho_k(H))$.
Then
\[
 |H:H^2| \> \leq \> 4^{s^2}.
\]
\end{lemma}

\begin{proof}
Let $P$ be a Sylow $2$-subgroup of $H$. The quotient $H/H^2$ is an
elementary abelian $2$-group, and the natural map $P \to H/H^2$ is surjective.
It follows that $|H:H^2| \leq 2^{d(P)}$.
Choose an elementary abelian subgroup $P_0 \leqslant P$ of maximal possible rank $r$.
By Corollary~\ref{cor:involutions} we obtain $\rho_k(H)\leq 2^{-r}$, which is precisely $r\leq s$.
A theorem of Mann \cite[Theorem~B]{Man71}
(see also \cite[Theorem~1.1(3) and Remark~1.2]{HPPS24}) gives
\[
 d(P)\leq r^2+\frac{r(r+1)}2\leq 2s^2.
\]
Combining the estimates gives the desired inequality.
\end{proof}

We are ready to prove the main results of this section.

\begin{proof}[Proof of Proposition~\ref{prop:abelianCreation}]
Let $H_0=H^2Z(H)$, and let $A=(H_0,\oplus)$ be the abelian group from
Lemma~\ref{lem:law}.
By Lemma~\ref{lem:coset} we have
\[ \rho_k(H) \leq \rho_k(H_0) = \rho_k(A) , \]
while Lemma~\ref{lem:index} gives $|A| \geq |H^2| \geq |H|\,4^{-s^2}$.
\end{proof}

\begin{proof}[Proof of Theorem~\ref{thm:abeTOclass2}]
Let $H$ be a finite nilpotent group of class $2$, and write
\[
 L=\log|H|,\qquad s=-\log_2(\rho_k(H)).
\]
If $s > \tfrac{1}{2}\sqrt{L}$, then
\[
 \rho_k(H)=2^{-s}
 \leq 2^{-\tfrac{1}{2}\sqrt{L}} ,
\]
and we are done because $\alpha \leq 1/2$.
Therefore we can assume $s \leq \tfrac{1}{2}\sqrt{L}$.

Let $A$ be the abelian group from Proposition~\ref{prop:abelianCreation},
so that
\[
 \log|A| \geq L-2(\log2)s^2 \geq \frac{L}{2} .
\]
If $c,C>0$ are constants witnessing $(\mathrm{AP}_{k,\alpha})$ for
finite abelian groups, then
\begin{align*}
 \rho_k(H)
 &\leq \rho_k(A)\\
 &\leq C\exp\bigl(-c(\log|A|)^\alpha\bigr)\\
 &\leq C\exp(-c \,2^{-\alpha} L^\alpha)
\end{align*}
as desired.
\end{proof}

Perhaps it is worth noting the bizarre coincidence that the proof of Theorem~\ref{thm:abeTOclass2} works precisely up to Behrend's threshold.

\vspace{0.1cm}
\section{Arbitrary groups}\label{sec:arbitrary}

Let $G$ be a finite group, and suppose the class of finite abelian groups satisfies $(\mathrm{AP}_{k,\alpha})$ for some $k \geq 3$ and $0<\alpha\leq 1/2$.
For constants $c,C,c_0>0$ and the subgroup $H \leqslant G$ supplied by \cite{Sab22} we have
\begin{align*}
 \rho_k(G)
 &\underset{\text{Lemma~\ref{lem:coset}}}{\leq}
 \rho_k(H)\\[1ex]
 &\underset{\text{Theorem~\ref{thm:abeTOclass2}}}{\leq}
 C\exp\bigl(-c(\log|H|)^\alpha\bigr)\\[1ex]
 &\underset{\text{Theorem~\ref{thm:large-subgroups}}}{\leq}
 C\exp\left(
 -c\,c_0^\alpha
 \left(\frac{\log|G|}{\log\log|G|}\right)^\alpha
 \right).
\end{align*}

The double logarithm is a nuisance.
In the rest of this section we explain the clean reduction to finite nilpotent groups, which is the most group-theoretical part of the paper.
Following \cite{Pyb97}, a finite group is of \emph{alternating type}
if all its nonabelian composition factors are alternating groups of
degree at least $8$ (smaller degrees are excluded for technical reasons). Note that solvable groups are included in this definition.
The {\em solvable radical} of a finite group is the largest solvable normal subgroup.

\begin{lemma}\label{lem:alternating-type}
Let $Q$ be a finite group of alternating type and trivial solvable radical.
Then $Q$ has an elementary abelian $3$-subgroup of order at least
$|Q|^{c/\log\log |Q|}$,
where $c>0$ is an absolute constant.
\end{lemma}

\begin{proof}
This was proved along \cite[Proof of Theorem~1]{Sab22}, but we rewrite a proof here for completeness.
Since the solvable radical of $Q$ is trivial, its socle has the form
\[
 D=\Soc(Q)=\prod_i\Alt(a_i)^{b_i},
 \qquad a_i\ge 8,
\]
where the $a_i$ are distinct.
Set $m=\sum_i a_i b_i$.
Disjoint $3$-cycles in each alternating factor give an elementary abelian $3$-subgroup of $D$ of cardinality at least $3^{m/4}$.

We have $C_Q(D)=1$, so $Q$ embeds into $\Aut(D)$ and
\[
 |D|\leq |Q|\leq\prod_i(a_i!)^{b_i}b_i! .
\]
Using $a_i,b_i\leq m$ in the upper bound and
$\log(a_i!/2)\gg a_i$ in the lower bound, we obtain
\[
 m \ll\log |Q|\ll m\log(2+m),
\]
and it follows that $m\gg\log |Q|/\log\log |Q|$.
\end{proof}

A {\em section} of $G$ is a quotient $M/N$ with $N\unlhd M\leqslant G$.
We want to use sections to estimate $r_k(G)$, so we first use a result of Gonz\'alez-S\'anchez and Klopsch \cite{GSK11} to obtain genuine subgroups on which to apply Lemma~\ref{lem:coset}.
In the next result, we choose the prime $p=3$ because \cite{GSK11} is a little weaker for $p=2$, and because in a group of exponent $3$ a progression repeats itself after the first three terms, so $r_k=r_3$ for any $k \geq 3$ in this case.

\begin{lemma}\label{lem:three-section}
There is an absolute constant $c>0$ such that, if a finite group
$G$ has a section $E$ which is an elementary abelian $3$-group, then
\[
 \rho_k(G) \> \ll \> |E|^{-c}
\]
for any $k \geq 3$.
\end{lemma}

\begin{proof}
Write $E=M/N\cong C_3^r$, where $N\unlhd M\leqslant G$, and choose a
Sylow $3$-subgroup $P$ of $M$. The projection of $P$ on $M/N$ is surjective, so
$d(P)\ge r$. Now \cite[Proposition~3.8]{GSK11} gives a subgroup $H\leqslant P$ of nilpotency class at most $2$ and exponent $3$ such that
\[
 |H|\ge |P:P^3[P,P]|=3^{d(P)}\ge3^r.
\]
By Lemma~\ref{lem:coset} and the observation above, we obtain
\[
 \rho_k(G) \leq \rho_k(H) = \rho_3(H) .
\]
Let $|H|=3^s$. Since $H$ has odd order and exponent $3$, the abelian group $(H=H_0,\oplus)$ supplied by the Baer trick is $\F_3^s$.
The cap-set bound \cite{EG17} gives
$r_3(\F_3^s) \ll \gamma^s$ for an absolute $\gamma<3$, and
using Lemma~\ref{lem:law} we obtain
\[
 \rho_3(H) = \frac{r_3(\F_3^s)}{3^s} \ll (\gamma/3)^s \leq |E|^{-c} 
\]
for some $c>0$, as desired.
\end{proof}

\begin{proof}[Proof of Theorem~\ref{thm:nilpTOarbitrary}]
By \cite[Corollary~2.3]{Pyb97}, there is an alternating-type subgroup
$T\leqslant G$ with $|T|\ge |G|^{1/5}$.
By Lemma~\ref{lem:coset} we have $\rho_k(G) \leq \rho_k(T)$, and we are going to bound the second.

Let $R$ be the solvable radical of $T$.
If $|R|\ge|T|^{1/2}$, then a well known consequence of P\'alfy--Wolf theorem \cite{Pal82,Wol82} supplies a nilpotent normal subgroup $H \unlhd R$ such that
\begin{equation} \label{eqNilpBig}
    |H| \geq |R|^{1/4} \geq |G|^{1/40} .
\end{equation}
If $c,C>0$ are constants witnessing $(\mathrm{AP}_{k,\alpha})$ for the class of finite
nilpotent groups, another application of Lemma~\ref{lem:coset} gives
\begin{align*}
 \rho_k(T)
 &\leq\rho_k(H)\leq C\exp\bigl(-c(\log|H|)^\alpha\bigr)\\
 &\leq C\exp\bigl(-c\,40^{-\alpha}(\log |G|)^\alpha\bigr) ,
\end{align*}
and we are done in this case.

Suppose otherwise that $|R|<|T|^{1/2}$, so that $|T/R| \ge |G|^{1/10}$.
Since $T/R$ has trivial solvable radical, by Lemma~\ref{lem:alternating-type} it has an elementary abelian $3$-subgroup of order at least $|G|^{c/\log\log|G|}$ for some constant $c>0$.
Up to replacing the constant $c$, by Lemma~\ref{lem:three-section} we have
\[
 \rho_k(T)
 \ll \exp\left(-c \frac{\log |G|}{\log\log |G|}\right) ,
\]
which is better than required.
This proves the first part.

If Conjecture~\ref{conj:A} holds for nilpotent groups with
exponent $\theta_k$, then (\ref{eqNilpBig}) gives
$\rho_k(G)\ll_k \rho_k(|G|)^{\theta_k/40}$.
In the second case, the conclusion follows from Behrend's
lower bound
$\rho_k(N)\ge\rho_3(N)\gg\exp(-c\sqrt{\log N})$.
\end{proof}

We point out that a positive answer to a question of Pyber \cite{Pyb97} (see \cite[Conjecture~6]{Sab22})
would close the gap between nilpotent groups of class $2$ and of arbitrary class, and further reduce Conjecture~\ref{conj:A} to the abelian case.

\vspace{0.1cm}
\section{Unconditional bounds for $k=3$} \label{sec:unconditional}

\subsection{From intervals to abelian groups} 

We now prove Conjecture~\ref{conj:A} for $k=3$ and abelian groups.
The following estimate is produced by embedding an arbitrary abelian group into a sufficiently large interval.
This argument is standard, see \cite[Proposition~5.24 and Lemma~5.25]{TV06}, and the novelty relies in combining it with \cite{EG17}.

\begin{proposition} \label{prop:intervTOabe}
There exists an absolute constant $\theta >0$ such that
for every finite abelian group $A$ we have
\[
    \rho_3(A) \> \ll \> \rho_3(|A|)^\theta .
\]
\end{proposition}

\begin{proof}
Let $N=|A|$, $r=d(A)$, and write
\[
 A\cong C_{m_1}\times\cdots\times C_{m_r},
 \qquad 2\leq m_1\mid m_2\mid\cdots\mid m_r .
\]

Set $w_1=1$, $w_{i+1}= 2 w_i m_i$ for every $1 \leq i \leq r-1$, and note that $w_r m_r = 2^{r-1}N$.
We define a function $\phi \colon A \to \{0,\ldots,2^rN-1\}$ by
\[ \phi(x_1,\ldots,x_r) :=
\sum_{i=1}^r w_i x_i = x_1 + 2m_1x_2 + 4m_1m_2x_3 + \cdots . \]
It is easy to check that $\phi$ is injective.
Now suppose $\phi(x)+\phi(z)=2\phi(y)$, which means $\sum_i w_i (x_i+z_i-2y_i) =0$.
Clearly the coefficients $c_i=x_i+z_i-2y_i$ satisfy $|c_i|<2m_i$.
Reducing modulo $2m_1$ gives $c_1=0$, and division by $2m_1$ followed by iteration gives $c_i=0$ for every $i$.
Therefore if $X\subseteq A$ is $3$-progression-free, then $\phi(X) \subseteq \{0,\ldots,2^rN-1\}$ is $3$-progression-free.

Partitioning $\{0,\ldots,2^rN-1\}$ into $2^r$ consecutive
intervals of length $N$ yields $|X| = |\phi(X)| \leq 2^r r_3(N)$.
Observe that $A$ contains an elementary abelian subgroup isomorphic to
$\F_p^r$ for some prime $p$.
By Lemma~\ref{lem:coset} and Theorem~\ref{thm:EG} we obtain
\[
 \rho_3(A)
 \leq \rho_3(\F_p^r)
 \ll e^{-c \, r} .
\]
It follows that
$2^r \leq C (N/|X|)^C$ for an absolute constant $C>0$.
Hence
\[
 \frac{|X|}{N} \leq C (N/|X|)^C \rho_3(N) ,
\]
and rearranging gives the claimed bound with $\theta=(1+C)^{-1}$.
\end{proof}

We now apply Proposition \ref{prop:intervTOabe} to the state-of-art result of Raghavan, and obtain our main unconditional estimate.

\begin{corollary}[Raghavan's bound for finite abelian groups] \label{cor:Raghavan}
There are absolute constants $c,C>0$ such that every finite abelian
group $A$ satisfies
\[
 r_3(A) \> \leq \> C |A| \exp\left(
 -c\left(\frac{\log\abs A}{\log\log\abs A}\right)^{1/6}
 \right).
\]
\end{corollary}
\begin{proof}
Combine \cite[Theorem~1.4]{Rag26} with Proposition \ref{prop:intervTOabe}.
\end{proof}

\begin{proof}[Proof of Theorem~\ref{thm:unconditional}]
The proof of Theorem~\ref{thm:abeTOclass2} applies equally to the bound in Corollary~\ref{cor:Raghavan}, so that estimate holds for finite nilpotent groups of class $2$ as well.
For the subgroup $H\leqslant G$ of Theorem~\ref{thm:large-subgroups} we have
\[ \frac{\log|H|}{\log\log|H|} \geq 
\frac{\log|H|}{\log\log|G|} \geq
c_0 \frac{\log|G|}{(\log\log|G|)^2} , \]
so by applying Lemma~\ref{lem:coset} we obtain
\begin{align*}
 \rho_3(G)
 &\leq \rho_3(H)\\
 &\leq C \exp\left(
 -c\left(\frac{\log\abs H}{\log\log\abs H}\right)^{1/6} \right) \\
 &\leq  C \exp \left( -c \, c_0^{1/6} \frac{(\log|G|)^{1/6}}{(\log\log|G|)^{1/3}} \right) .
\end{align*}
The proof follows.
\end{proof}

\subsection{Quasirandom groups}

Quasirandomness, as introduced by Gowers \cite{Gow08}, gives bounds for $3$-progression-free sets that can be
substantially stronger than Theorem~\ref{thm:unconditional}.

For example, Tao \cite{Tao13} proved mixing estimates for progressions
in $\mathrm{SL}_d(\mathbb F_q)$ with $d$ fixed.
Peluse \cite{Pel18} obtained similar estimates for all
nonabelian finite simple groups,
and Bhangale, Harsha and Roy \cite{BHR22} extended the mixing result to
abstract quasirandom groups. Their theorem has the following
immediate consequence.

\begin{proposition}[Bhangale--Harsha--Roy]
Let $G$ be a finite nontrivial group and suppose that every nontrivial irreducible
complex representation of $G$ has dimension at least $D \geq 1$.
Then
\[
 r_3(G) \> \ll \> \abs G\,D^{-1/24}.
\]
\end{proposition}

\begin{proof}
Let $X\subseteq G$ be $3$-progression-free. By \cite[Theorem~1]{BHR22},
\[
 \frac{|X|^3}{|G|^3} \leq\E_{x,y\in G}
 [\one_X(x)\one_X(xy)\one_X(xy^2)]
       +2^{1/4}D^{-1/8} . \]
       Since only $y=1$ contributes to the expectation, we obtain
 \[ \frac{|X|^3}{|G|^3} \leq \frac{|X|}{|G|} \cdot \frac{1}{|G|} +2^{1/4}D^{-1/8} \leq \frac{1}{|G|} +2^{1/4}D^{-1/8} . \]
The sum-of-squares formula for irreducible character degrees gives
$|G| > D^2$. It follows that
\[
 \frac{|X|^3}{|G|^3} < D^{-2}+2^{1/4}D^{-1/8}\ll D^{-1/8}
\]
as desired.
\end{proof}

We refer the reader to \cite{BS25} for an algebraic description of highly quasirandom groups;
in particular, if $D \geq c \log|G|$, then the quotient of $G$ by its solvable radical is a direct product of nonabelian finite simple groups \cite[Remark~4.6]{BS25}.

\vspace{0.2cm}

\end{document}